\documentclass[12pt, a4paper]{article}
\usepackage[ansinew]{inputenc}
\usepackage[T1]{fontenc} 
\usepackage{longtable}
\usepackage{xcolor}
\usepackage{graphicx}
\usepackage{setspace}
\usepackage{amssymb,amsmath}
\usepackage[all]{xy}
\usepackage{lscape}
\usepackage{hyperref}

\newcounter{StructureNumber}
\newenvironment{Structure}[1]{\noindent \bf #1 \arabic{StructureNumber}. \rm}{\stepcounter{StructureNumber}}
\newenvironment{Proof}{\noindent \bf Proof. \rm}{\hfill $\diamondsuit$}

\begin{document}
\title{
On Dirichlet Products Evaluated at Fibonacci Numbers
}

\author{Uwe Stroinski\footnote{E-mail address: uwst@fa.uni-tuebingen.de. 
I thank Rainer Nagel for his valuable sup\-port and the Tuebingen Functional Analysis group for their kind hospitality.}}

\maketitle

\onehalfspacing

\begin{abstract}
In this work we discuss sums of the form $\sum_{n\leq x} (f*g)(a_n)$ with $a_n$ being the $n$--th Fibonacci number.
As first applications of the results we get a representation of Fibonacci numbers in terms of Euler's $\varphi$--function, an upper bound on the number of primitive prime divisors and a non--trivial fixed point $h(m) = \sum_{\alpha(n)= m} h(n)$, where $\alpha(n)$ is the least index $m$ such that $n|a_m$.
\end{abstract}

 \section{Introduction}

Let $a_n$ denote the Fibonacci sequence, thus $a_n$ is recursively given by $a_0:=0,a_1:=1$ and $a_{n+1}:=a_n+a_{n-1}$ for $n\geq1$.
In the following we consider a version of the well--known formula from multiplicative number theory
\[
\sum_{n\leq x} (f*g) (n) = \sum_{n\leq x} f(n) G\left(\frac{x}{n}\right)= \sum_{n\leq x} g(n) F\left(\frac{x}{n}\right)
\]
where $G(x) := \sum_{n\leq x} g(n)$ and  $F(x) := \sum_{n\leq x} f(n)$ adapted to Dirichlet products evaluated at Fibonacci numbers.
\\
\pagebreak

\begin{Structure}{Theorem}
\newcounter{DoubleCounting}
\setcounter{DoubleCounting}{\value{StructureNumber}}
Let $f$ and $g$ be arithmetic functions and let $x\geq 1$ be real.
Then we have
\[
\sum_{n\leq x} (f*g)(a_n)
=
\sum_{\alpha(n)\leq x} f(n) \sum_{d\leq \frac{x}{\alpha(n)}} g\left(\frac{a_{d\cdot \alpha(n)}}{n}\right)
=
\sum_{\alpha(n)\leq x} g(n) \sum_{d\leq \frac{x}{\alpha(n)}} f\left(\frac{a_{d\cdot \alpha(n)}}{n}\right).
\]
\end{Structure}
\\

\noindent
Here $\alpha$ denotes the dual of the Fibonacci sequence (sometimes called the rank of apparition) and is defined as 
\[
\alpha(n):=\min\{m\in\mathbb{N}: n|a_m\},
\]
the smallest integer $m$ such that $n|a_m$.
\\

\noindent
I discuss the relation of this theorem to known results of Matiyasevich and Guy \cite{MatGuy86} as well as results of Kiss \cite{Kiss88} and give some further applications.
Among them are a new representation of Fibonacci numbers in terms of Euler's $\varphi$--function, an upper bound for the number of primitive prime factors in the Fibonacci sequence and explicit representations of some related Euler products.
\\

\noindent
Let me fix some notation.
Define
\[
e_n:=\max\{m\in\mathbb{N}:n^m|a_{\alpha(n)}\},
\]
the maximal power of $n$ that still divides $a_{\alpha(n)}$.
The integer part of a real number $x$ is $\lfloor x\rfloor$ and the fractional part is $\{x\}$.
Therefore $x=\lfloor x\rfloor +\{x\}$.  
Define the golden ratio $r:=\frac{1+\sqrt{5}}{2}$ as well as $s:=\frac{1-\sqrt{5}}{2}$ and remember Binet's formula $a_n=\frac{r^n-s^n}{\sqrt{5}}$. 
An arithmetic function is a function from $\mathbb{N}$ to $\mathbb{C}$.
The sum $\sum_{\alpha(n)\leq x}f(n)$ means $\sum_{n:{1\leq n \atop\alpha(n)\leq x}}f(n)$, that is we sum over all $n \geq 1$ with $\alpha(n)\leq x$.
These sums are finite since for $n>a_{\lfloor x\rfloor}$ one has $\alpha(n)>\lfloor x\rfloor$.
In what follows $p$ is always a prime number and $\sum_{\alpha(p)\leq x}f(p)$ denotes the sum of $f(p)$ over all primes that divide at least one Fibonacci number $a_n$ with $n\leq x$. 
\\

\noindent
For an arithmetic function $f$ define its $\alpha$--contraction $f_\alpha$ as
\[
f_\alpha(n):=\sum_{m:\alpha(m)=n}f(m)
\]
with the empty sum defined to be zero $f_\alpha(2):=0$. 
Set 
\[
T_{f,\alpha}(x):=\sum_{\alpha(n)\leq x}f(n)\left\lfloor\frac{x}{\alpha(n)}\right\rfloor
\]
and 
\[
S_{f,\alpha}(x):=\sum_{\alpha(n)\leq x} f(n).
\]
Observe that this implies $S_{f,\alpha}(x)=\sum_{n\leq x}f_\alpha(n)$.
Remember, that a standard result now yields
$
T_{f,\alpha}(x)
=
\sum_{n\leq x}S_{f,\alpha}\left(\frac{x}{n}\right)
$
and by Moebius inversion we obtain
$
S_{f,\alpha}(x)
=
\sum_{n\leq x}\mu(n)T_ {f,\alpha}\left(\frac{x}{n}\right).
$
We keep the $\alpha$ in the index since most of the results hold with exactly the same proofs for more general strong divisibility sequences.
Since we focus on Fibonacci numbers, we do not feel the need to present each preliminary result in most general form.
Especially since it is Fibonacci numbers that have recently received attention (eg. a polymath proposal concerning Littlewood's conjecture by Gowers \cite{Gow09}) and since the problems are already considered hard for this special case (cf. eg. Sarnak \cite{Sar07} on the number of Fibonacci primes or Granville \cite{Gran12} on the existence of Wall--Sun--Sun primes).

\section{The proof of Theorem \arabic{DoubleCounting}}

\noindent
In the next lemma we collect some results on the dual sequence.
\\

\begin{Structure}{Lemma}
Let $0<n,m\in\mathbb{N}$ and $p$ be prime, then we have
\begin{itemize}
\item duality: $n|a_m$ if and only if $\alpha(n)|m$,
\item $\alpha(p^n) = p^{n-e_p}\alpha(p)$ for odd primes $p$ with $n>e_p$ and
\item $\alpha(2)=3, \alpha(4)=6$ and $\alpha(2^n) = 3 \cdot 2^{n-2}$ for $n\geq 3$.
\end{itemize}
\end{Structure}

\noindent
The proof of this lemma is well--known and can be found e.g in Halton \cite{Hal66}.
\\

\noindent
Let us now prove Theorem \arabic{DoubleCounting}.
\\

\begin{Proof}
Let $1\leq x \in\mathbb{R}$ and $n=\lfloor x\rfloor$.
For each $x$ we consider a table with $a_n$ rows and $n$ columns.
Let $1\leq d \leq a_n$ and fill position $(d,n)$ with $f(d)g\left(\frac{a_n}{d}\right)$ if $\alpha(d)|n$ and with zero otherwise. 
This table is depicted in Figure 1.
\\

\noindent
Double counting yields
\[
\sum_{n\leq x} \sum_{d:\alpha(d)|n}f(d) g\left(\frac{a_n}{d}\right)= \sum_{n\leq a_x}f(n)\sum _{d \leq x/\alpha(n)}  g\left(\frac{a_{d\cdot \alpha(n)}}{n}\right).
\]
For $\alpha(n)> x$ we have $\sum _{d \leq x/\alpha(n)}  g\left(\frac{a_{d\cdot \alpha(n)}}{n}\right)=0$ and therefore the summation in the right sum can be restricted to those $n\leq a_x$ with $\alpha(n)\leq x$ and since $n>a_x$ implies $\alpha(n)>x$ it suffices to sum over those $n$ with $\alpha(n)\leq x$.
In the left sum we observe that $d|a_n$ if and only if $ \alpha(d)|n$.
That proves the first equality.
The second equality in the theorem follows from $f*g = g*f$.
\end{Proof}
\\

\noindent
The applicability of Theorem \arabic{DoubleCounting} depends on the ability to control sums like $\sum_{d\leq \frac{x}{\alpha(n)}} g\left(\frac{a_{d\cdot \alpha(n)}}{n}\right)$.
In general, this is hard.  
However, one frequently considers functions $g$ satisfying additional assumptions.
\\

\begin{Structure}{Corollary}
\newcounter{DirichletProdAtFibo}
\setcounter{DirichletProdAtFibo}{\value{StructureNumber}}
Let $f$ and $g$ be arithmetic functions and let $g$ be completely multiplicative.
Then we have
\[
(f*g)(a_n)
=
g(a_n) \left(1*\left(\frac{f}{g}\right)_\alpha\right)(n) .
\]
If additionally $g(n) =1$ for all $n\in\mathbb{N}$ we have
\[
(1*f)(a_n) = (1*f_\alpha)(n).
\]
\end{Structure}

\begin{landscape}
\begin{figure}[h]
\[
\begin{array}{c||c|c|c|c|c|c|c|c||c }
\hline
n & 1 & 2 & 3 & 4 & 5 & 6 & \ldots & n &\\
\hline
d \backslash a_n & 1 & 1 & 2 & 3 & 5 & 8 &\ldots & a_n & \sum\\
\hline\hline
1 & f(1) g\left(\frac{1}{1}\right) & f(1) g\left(\frac{1}{1}\right) & f(1) g\left(\frac{2}{1}\right) & f(1) g\left(\frac{3}{1}\right) & f(1) g\left(\frac{5}{1}\right) & f(1) g\left(\frac{8}{1}\right) & \cdots & f(1 )g\left(\frac{a_n}{1}\right) & f(1) \sum _{d \leq n/\alpha(1)}  g\left(\frac{a_{d\cdot \alpha(1)}}{1}\right) \rule{0mm}{0.5cm}\\[0.5ex]
\hline
2 & 0 & 0 & f(2) g\left(\frac{2}{2}\right) & 0 &  0 & f(2) g\left(\frac{8}{2}\right) & \cdots & \cdots & f(2)\sum _{d \leq n/\alpha(2)}  g\left(\frac{a_{d\cdot \alpha(2)}}{2}\right) \rule{0mm}{0.5cm}\\[0.5ex]
\hline
3 & 0 & 0 & 0 & f(3) g\left(\frac{3}{3}\right) & 0 & 0 & \cdots & \cdots & f(3)\sum _{d \leq n/\alpha(3)}  g\left(\frac{a_{d\cdot \alpha(3)}}{3}\right) \rule{0mm}{0.5cm}\\[0.5ex]
\hline
4 & 0 & 0 & 0 & 0 & 0 & f(4) g\left(\frac{8}{4}\right) & \cdots & \cdots & f(4)\sum _{d \leq n/\alpha(4)}  g\left(\frac{a_{d\cdot \alpha(4)}}{4}\right) \rule{0mm}{0.5cm}\\[0.5ex]
\hline
5 & 0 & 0 & 0 & 0 & f(5) g\left(\frac{5}{5}\right)& 0 & \ldots & \ldots  & f(5)\sum _{d \leq n/\alpha(5)}  g\left(\frac{a_{d\cdot \alpha(5)}}{5}\right) \rule{0mm}{0.5cm}\\[0.5ex]
\hline
\vdots & \vdots & \vdots & \vdots & \vdots & \vdots & \vdots & \vdots & \vdots & \vdots \rule{0mm}{0.5cm} \\[0.5ex]
\hline
a_n & 0 & 0 & 0 & 0 & 0 & 0 &\ldots & f(a_n) g\left(\frac{a_n}{a_n}\right) & f(a_n)\sum _{d \leq n/\alpha(a_n)}  g\left(\frac{a_{d\cdot \alpha(a_n)}}{a_n}\right)\rule{0mm}{0.5cm}\\[0.5ex]
\hline
\hline
\sum & \displaystyle \sum_{\alpha(d)|1}f(1)g\left(\frac{1}{d}\right) & \ldots & \ldots & \ldots & \ldots & \ldots &\ldots & \displaystyle \sum_{\alpha(d)|n}f(d)g\left(\frac{a_n}{d}\right) \rule{0mm}{0.75cm} & \\
\hline
\end{array}
\]
\caption{The double counting argument}
\end{figure}
\end{landscape}

\begin{Proof}
By Theorem \arabic{DoubleCounting} we have
\begin{eqnarray*}
\sum_{n\leq x} (f*g)(a_n)
& = &
\sum_{\alpha(n)\leq x} f(n) \sum_{d\leq \frac{x}{\alpha(n)}} g\left(\frac{a_{d\cdot \alpha(n)}}{n}\right).
\end{eqnarray*}
Since $g$ is completely multiplicative, this equals
\begin{eqnarray*}
& = &
\sum_{d,n: d\cdot \alpha(n)\leq x} \frac{f(n)}{g(n)}  g(a_{d\cdot \alpha(n)}).
\end{eqnarray*}
We rearrange the sum
\begin{eqnarray*}
& = &
\sum_{d, k : d\cdot k\leq x} g(a_{d\cdot k}) \sum_{n:\alpha(n)=k} \frac{f(n)}{g(n)}
\end{eqnarray*}
and use the definition of $\alpha$--contraction
\begin{eqnarray*}
& = &
\sum_{d, k : d\cdot k\leq x} g(a_{d\cdot k}) \left(\frac{f}{g}\right)_\alpha(k).
\end{eqnarray*}
Introducing a new index $n= d\cdot k$ yields
\begin{eqnarray*}
& = &
\sum_{n\leq x} g(a_n) \sum_{k|n} \left(\frac{f}{g}\right)_\alpha(k)
\end{eqnarray*}
and, by the definition of the Dirichlet product, we get
\begin{eqnarray*}
\sum_{n\leq x} (f*g)(a_n)
& = &
\sum_{n\leq x} g(a_n) \left(1* \left(\frac{f}{g}\right)_\alpha\right)(n).
\end{eqnarray*}
These sums are equal for all $x\in\mathbb{R}$ and therefore
\[
(f*g)(a_n) = g(a_n) \left(1* \left(\frac{f}{g}\right)_\alpha\right)(n)
\]
for $n\in\mathbb{N}$. 
That proves the first assertion.
The second assertion is an immediate consequence. 
\end{Proof}

\section{Relation to existing results}

\noindent
In this section we show how the above connects to existing work (e.g. Matiyasevich and Guy \cite{MatGuy86} or Kiss \cite{Kiss88}) and prepare the next steps. 
\\

\begin{Structure}{Corollary}
\newcounter{LambdaSum}
\addtocounter{LambdaSum}{\value{StructureNumber}}
For $x\geq 1$ we have
\begin{eqnarray*}
\sum_{\alpha(n)\leq x}\Lambda(n)\left\lfloor\frac{x}{\alpha(n)}\right\rfloor
& = &
\frac{\log r}{2}\lfloor x\rfloor^2
+
\frac{\log \frac{r}{5}}{2}\lfloor x\rfloor
+
\sum_{n\leq x} \log \left(1-\frac{(-1)^n}{r^{2n}}\right).
\end{eqnarray*}
\end{Structure}
\noindent
Note that
\[
\sum_{\alpha(n)\leq x}\Lambda(n)\left\lfloor\frac{x}{\alpha(n)}\right\rfloor = \log \prod_{n\leq x} a_n.
\]

\begin{Proof}
By Theorem \arabic{DoubleCounting} one has 
$\sum_{\alpha(n)\leq x}\Lambda(n)\left\lfloor\frac{x}{\alpha(n)}\right\rfloor=\sum_{n\leq x} \sum_{d|a_n}\Lambda(d)$
and since 
$\sum_{d|n}\Lambda(d)=\log n$ this yields $\sum_{\alpha(n)\leq x}\Lambda(n)\left\lfloor\frac{x}{\alpha(n)}\right\rfloor=\log \prod_{n\leq x} a_n$.
Binet's formula together with $s r = -1$ implies
\begin{eqnarray*}
\log \left[\prod_{n\leq x} a_n\right]
& = &
\log \left[\left(\frac{1}{\sqrt{5}}\right)^{\lfloor x\rfloor}\prod_{n\leq x}r^n\prod_{n\leq x}1-\frac{s^n}{r^n}\right]\\
& = &
\log \left[\left(\frac{1}{\sqrt{5}}\right)^{\lfloor x\rfloor}\prod_{n\leq x}r^n\prod_{n\leq x}1-\frac{(-1)^n}{r^{2n}}\right]\\
& = &
\frac{\log r}{2}\lfloor x\rfloor^2+\frac{\log \frac{r}{5}}{2}\lfloor x\rfloor+\sum_{n\leq x} \log \left(1-\frac{(-1)^n}{r^{2n}}\right)
\end{eqnarray*}
which proves the corollary.
\end{Proof}
\\

\noindent
Observe that 
\[
c:=\lim_{x\rightarrow\infty}\sum_{n\leq x} \log \left(1-\frac{(-1)^n}{r^{2n}}\right)
\]
exists with approximately $c=0.2043618834...$.
\\

\noindent
The next result delivers the exact asymptotic growth of the summatory (in the sense of this paper) von Mangoldt function.
\\

\begin{Structure}{Corollary}
\newcounter{AsympMangoldt}
\setcounter{AsympMangoldt}{\value{StructureNumber}}
We have
\[
\sum_{\alpha(n)\leq x}\Lambda(n) = \frac{3 \log r}{\pi^2}x^2+O(x\log x).
\]
\end{Structure}
\noindent
Note that
\[
\sum_{\alpha(n)\leq x}\Lambda(n) = \log \textnormal{lcm}(a_1,a_2,\ldots, a_{\lfloor x\rfloor}).
\]
\\

\begin{Proof}
With the notations of the previous section we obtain for $f=\Lambda$ and $c_x:=\sum_{n\leq x} \log \left(1-\frac{(-1)^n}{r^{2n}}\right)$ that
\begin{eqnarray*}
S_{\Lambda,\alpha}(x)
& = &
\sum_{n\leq x}\mu(n)T_ {\Lambda,\alpha}\left(\frac{x}{n}\right)\\
& = &
\frac{\log r}{2}\sum_{n\leq x}\mu(n)\left\lfloor \frac{x}{n}\right\rfloor^2
+
\frac{\log \frac{r}{5}}{2}\sum_{n\leq x}\mu(n)\left\lfloor \frac{x}{n}\right\rfloor
+
c_x\sum_{n\leq x}\mu(n) \\
& = &
\frac{\log r}{2}\sum_{n\leq x}\mu(n)\left\lfloor \frac{x}{n}\right\rfloor^2
+
c_x\sum_{n\leq x}\mu(n)
+
\frac{\log \frac{r}{5}}{2}.
\end{eqnarray*}
The sum $\sum_{n\leq x}\mu(n)$ trivially is $O(x)$ and it is an exercise to show that 
\[
\sum_{n\leq x}\mu(n)\left\lfloor \frac{x}{n}\right\rfloor^2 = \frac{1}{\zeta(2)}x^2 + O(x\log x).
\]
Since $\zeta(2) = \frac{\pi^2}{6}$ this finishes the proof.
\end{Proof}
\\

\noindent
By considering the quotient of the quantities in Corollary \arabic{LambdaSum} and Corollary \arabic{AsympMangoldt} and taking limits the formula for $\pi$ in \cite{MatGuy86} follows.
These ideas can  be generalized to get representations for $\zeta(k)$ for values $k$ not necessarily $2$ (cf. Akiyama \cite{Aki93}). 
\\

\noindent
Let us further investigate the summatory (in the sense of this paper) von Mangoldt function.
To that purpose we keep in mind that $2$ is the exception to $\alpha(p^n) = p^{n-e_p}\alpha(p)$ with $n>e_p$ and that there might exist primes with $e_p>1$
(cf. eg. Wall \cite{Wall60} as well as Sun and Sun \cite{SunSun92}).
We get
\begin{eqnarray*}
\sum_{\alpha(n)\leq x} \Lambda(n)
 & = &
\sum_{n=1}^\infty \sum_{\alpha(p^n)\leq x} \log p \\
 & = &
\sum_{\alpha(p)\leq x}\sum_{n=1}^{e_p} \log p 
+
\sum_{\substack{2 < p, e_p < n \\ p^{n-e_p}\alpha(p)\leq x}} \log p
+ 
\log 2\sum_{1 < n, \alpha(2^n)\leq x}1 
\\
\end{eqnarray*}
and conclude for $x\geq1$
\[
\sum_{\alpha(n)\leq x} \Lambda(n)
-
\sum_{\alpha(p)\leq x} e_p \log p 
=
\sum_{\substack{2 < p, e_p < n \\ p^{n-e_p}\alpha(p)\leq x}} \log p
+
O(\log x).
\]
Since in the right sum $n>e_p$ we have $p\leq x$ and thus
\[
\sum_{\substack{2 < p, e_p < n \\ p^{n-e_p}\alpha(p)\leq x}} \log p 
\leq
\sum_{p^n\leq x} \log p 
=
\sum_{n\leq x}\Lambda(n) 
\]
which is known to be $O(x)$.
Thus
\[
\sum_{\alpha(n)\leq x} \Lambda(n)
-
\sum_{\alpha(p)\leq x} e_p \log p 
= 
O(x)
\]
and we get a special case of a result due to Kiss \cite{Kiss88}.
\\

\begin{Structure}{Corollary}
\newcounter{AsympEP}
\setcounter{AsympEP}{\value{StructureNumber}}
We have
\[
\sum_{\alpha(p)\leq x}e_p\log p = \frac{3 \log r}{\pi^2} x^2 +O(x\log x).
\]
\end{Structure}

\section{An upper bound for the number of primitive primes}

\noindent
A primitive prime factor of $a_n$ is a prime $p$ that divides $a_n$ but does not divide $a_m$ for $1 \leq m \leq n - 1$.
The existence of primitive prime divisors for Lucas numbers with positive discriminant and suffciently large index (with Fibonacci numbers as a special case) is ensured by a result of Lucas \cite{Luc78} that was later considerably generalized (cf. Carmichael \cite{Carmichael13}, Schinzel \cite{Schinzel74}, Bilu, Hanrot and Voutier \cite{BilHanVou01}).
\\

\noindent 
Once existence is established it is natural to ask how many primitive prime factors there are.
We define
\[
\pi_\alpha(x):= \sum_{\alpha(p)\leq x}1
\]
as the number of primitive prime factors of Fibonacci numbers up to index $\lfloor x\rfloor$ and go on to show how Corollary \arabic{AsympMangoldt} can be used to get an upper bound complementing the existence results.
\\

\begin{Structure}{Theorem}
\newcounter{PiAsymp}
\setcounter{PiAsymp}{\value{StructureNumber}}
We have
\[
\limsup_{x\rightarrow\infty}\frac{\log x}{x^2}\pi_{\alpha}(x)\leq \frac{3 \log r}{2\pi^2}.
\]
\end{Structure}

\begin{Proof}
Fix $0<\varepsilon<2$ and write $\pi_\alpha$ as follows 
\begin{eqnarray*}
\pi_{\alpha}(x)
 & = &
\sum_{p \leq x^{2-\varepsilon}:\alpha(p)\leq x}1
+
 \sum_{p > x^{2-\varepsilon}:\alpha(p)\leq x}1.
\end{eqnarray*}
The first sum $\sum_{p \leq x^{2-\varepsilon}:\alpha(p)<x}1$ can be estimated by $x^{2-\varepsilon}$.
To get an estimate for the second sum let $\{p_1,\ldots,p_t\} := \{x^{2-\varepsilon} < p \leq a_{\lfloor x\rfloor}: \alpha(p)\leq x\}$.
Then $\prod_{n=1}^t p_n \leq \textnormal{lcm}_{n=1}^x a_n$.
Let $s$ be such that $\left(x^{2-\varepsilon}\right)^s = \textnormal{lcm}_{n=1}^x a_n$.
Since $x^{2-\varepsilon} < p_i$ for $1\leq i\leq t$ we get $t\leq s$ and thus
\[
\left(x^{2-\varepsilon}\right)^t
\leq 
\textnormal{lcm}_{n=1}^x a_n.
\]
Therefore
\[
\pi_{\alpha}(x)
\leq
x^{2-\varepsilon}+
\frac{\log \left(\textnormal{lcm}_{n=1}^x a_n\right)}{(2-\varepsilon)\log x}
\]
and with Corollary \arabic{AsympMangoldt} we get 
\[
\pi_{\alpha}(x)
\leq
x^{2-\varepsilon}
+
\frac{3 \log r}{(2-\varepsilon)\pi^2}\frac{x^2}{\log x}.
\]
This inequality is valid for any $0< \varepsilon < 2$ and therefore finished the proof of the theorem.
\end{Proof}
\\

\noindent
Let me give an alternative proof of this theorem.
\\

\begin{Proof}
Let $d(n)$ denote the number of divisors of $n$.
We have
\[
d\left(\prod_{\alpha(p)\leq x}p^{e(p)}\right)
=
\prod_{\alpha(p)\leq x} (1+e(p))
\]
and we know $\limsup_{n\rightarrow\infty}\frac{\log d(n) \log \log n}{\log n}=\log 2.$
Therefore Corollary \arabic{AsympMangoldt} implies
\[
\limsup_{x\rightarrow\infty}\frac{\sum_{\alpha(p)\leq x} \log(1+e(p)) 2 \log x}{\frac{3 \log r}{\pi^2}x^2}
\leq
\log 2
\]
and thus
\[
\limsup_{x\rightarrow\infty}\frac{\log x}{x^2}\sum_{\alpha(p)\leq x} \log(1+e(p))
\leq
\frac{3 \log r}{2 \pi^2}
\log 2.
\]
This implies
\[
\limsup_{x\rightarrow\infty}\frac{\log x}{x^2}\pi_\alpha(x)
\leq
\frac{3 \log r}{2 \pi^2}
\]
and thus completes the proof.
\end{Proof}
\\

\noindent
There is still a gap between the lower bounds indicated by Carmichael's theorem and the upper bound given above.
One might wonder why not to directly prove a prime number theorem giving the exact asymptotics.
The difficulty here stems from the fact that current proofs of the prime number theorem make use of $\sum_{n\leq x}\Lambda(n) \sim \sum_{p\leq x} \log p$ in a crucial way.
In our situation we have $\sum_{\alpha(n)\leq x} \Lambda(n) \sim \sum_{\alpha(p)\leq x} e_p \log p$ and $\sum_{\alpha(p) \leq x}\log \alpha(p) = \pi_\alpha(x)\log x - \int_1^x \frac{\pi_\alpha(t)}{t} dt$ by partial summation.
However, in the light of the above estimate, $|\sum_{\alpha(p)\leq x} e_p \log p-\sum_{\alpha(p) \leq x}\log \alpha(p)|= |\sum_{\alpha(p)\leq x}\log\frac{p^{e_p}}{\alpha(p^{e_p})}|$ is not asymptotically zero.
The connection of the growth of $\sum_{\alpha(n)\leq x} \Lambda(n)$ and the growth of $\pi_\alpha(x)$ is not as tight as in the situation of the prime number theorem.  

\section{A representation of Fibonacci numbers}

\noindent
Let $\varphi$ be Euler's totient function.
Since $\sum_{d|n}\varphi(d)=n$ we obtain from Theorem \arabic{DoubleCounting} that
\[
\label{phiSum}
\sum_{\alpha(n)\leq x}\varphi(n)\left\lfloor\frac{x}{\alpha(n)}\right\rfloor=\sum_{n\leq x}a_n = a_{\lfloor x\rfloor+2}-1
\]
where the last equality follows from an elementary induction argument.
\\

\noindent
Note that this can be used as an alternative definition of the Fibonacci numbers.
To that purpose set $a_1:=1$ and define 
\[
a_{x+2}:=1+\sum_{\alpha(n)\leq x}\varphi(n)\left\lfloor\frac{x}{\alpha(n)}\right\rfloor
\]
for $x\geq0$. 
In that case $\alpha$ is also defined recursively, the sum is taken over all numbers dividing at least one $a_n$ with $n\leq x$.
\\

\noindent
Euler's totient function appears in many places in number theory.
A connection to Fibonacci numbers like the above seems to be new.

\section{Iterating $\alpha$--contractions of Moebius $\mu$}

\noindent
In this section section we will repeatedly apply $\alpha$--contraction starting with the Moebius $\mu$ function.
As a result we will get a non--trivial fixed point and a non--trivial kernel element of $\alpha$--contraction considered as linear operator on arithmetic functions.
\\

\noindent
We know that $\sum_{d|n}\mu(d)=\left\lfloor\frac{1}{n}\right\rfloor$.
This expression is equal to $1$ if $n=1$ and equal to $0$ else. 
Since $a_1=a_2=1<a_n$ for $n>2$ we can evaluate $\sum_{n\leq x} \left\lfloor\frac{1}{a_n}\right\rfloor$ and if we now apply Theorem \arabic{DoubleCounting} we obtain
\begin{eqnarray*}
T_{\mu,\alpha}(x) 
=
\sum_{\alpha(n)\leq x}\mu(n)\left\lfloor\frac{x}{\alpha(n)}\right\rfloor
& = &
\left\{
\begin{array}{ll}
1, & 1 \leq x < 2,\\
2, & 2 \leq x.
\end{array}
\right.
\end{eqnarray*}
Moebius inversion yields
\begin{eqnarray*}
S_{\mu,\alpha}(x)
& = & 
\sum_{n\leq x}\mu(n)T_ {\mu,\alpha}\left(\frac{x}{n}\right) \\
& = &
2 \sum_{n\leq \frac{x}{2}} \mu(n) + \sum_{\frac{x}{2}<n\leq x}\mu(n) \\
& = &
\sum_{n\leq x}\mu(n) + \sum_{n\leq \frac{x}{2}}\mu(n)
\end{eqnarray*}
We set $M(x):=  \sum_{n\leq x}\mu(n)$ and finally arrive at
\[
S_{\mu,\alpha}(x)
=
M(x) + M\left(\frac{x}{2}\right).
\] 
We use that to get an explicit representation of $\mu_\alpha$.
\begin{eqnarray*}
\mu_\alpha(m)
& = & 
\sum_{\alpha(n)=m}\mu(n) \\
& = &
\sum_{\alpha(n)\leq m}\mu(n) - \sum_{\alpha(n)\leq m-1}\mu(n) \\
& = & 
S_{\mu,\alpha}(n) - S_{\mu,\alpha}(n-1) \\
& = & 
M(n) -M(n-1) + M\left(\frac{n}{2}\right) - M\left(\frac{n-1}{2}\right)\\
& = &
\mu(n) + 
\left\{
\begin{array}{ll}
\mu\left(\frac{n}{2}\right), & n \textnormal{ even},\\
0, & n \textnormal{ odd}.
\end{array}
\right.
\end{eqnarray*}

\noindent
Therefore
\[
\mu_\alpha(n) 
= 
\left\{
\begin{array}{ll}
\mu(n), & \textnormal{for } n\equiv 1,3 \mod 4\\
\mu\left(\frac{n}{2}\right), & \textnormal{for } n\equiv 0 \mod 4\\
0 , & \textnormal{for } n\equiv 2 \mod 4\\
\end{array}
\right.
\]
This function is multiplicative.
\\

\noindent
Define $\mu_{\alpha^{m+1}}:=(\mu_{\alpha^m})_\alpha$ for $m\geq 1$ and observe that
\begin{eqnarray*}
T_{\mu_\alpha,\alpha}(x)
& = &
\sum_{\alpha(n)\leq x} \mu_\alpha(n) \left\lfloor \frac{x}{\alpha(n)} \right\rfloor \\
& = &
\sum_{\alpha(n)\leq x} (1*\mu_\alpha)(a_n) \\
& = & 
\sum_{\alpha(n)\leq x} \left\lfloor \frac{1}{a_{a_n}} \right\rfloor \\
& = &
\left\{
\begin{array}{ll}
1, & n < 2, \\
2, & n < 3, \\
3, & n\geq 3.
\end{array}
\right.
\end{eqnarray*}
Therefore
\begin{eqnarray*}
S_{\mu_\alpha,\alpha}(x)
& = &
\sum_{n\leq x} \mu(n) T_{\mu_\alpha,\alpha}\left(\frac{x}{n}\right) \\
& = &
3 \sum_{n\leq x/3} \mu(n)  
+
2 \sum_{x/3 < n\leq x/2} \mu(n)
+
\sum_{x/2 < n\leq x} \mu(n)\\
& = &
M(x) + M\left(\frac{x}{2}\right) + M\left(\frac{x}{3}\right)
\end{eqnarray*}
and thus
\begin{eqnarray*}
\mu_{\alpha^2}(n)
& = &
M(n) - M(n-1) + M\left(\frac{x}{2}\right) - M\left(\frac{x-1}{2}\right) + M\left(\frac{x}{3}\right)- M\left(\frac{x-1}{3}\right) \\
& = &
\mu(n)
+
\left\{
\begin{array}{ll}
\mu\left(\frac{n}{2}\right), & n\equiv 0\mod 2, \\
0, & n\not\equiv 0\mod 2,
\end{array}
\right.
+
\left\{
\begin{array}{ll}
\mu\left(\frac{n}{3}\right), & n\equiv 0\mod 3, \\
0, & n\not\equiv 0\mod 3,
\end{array}
\right.\\
& = &
\mu(n)
+
\left\{
\begin{array}{ll}
\mu\left(\frac{n}{2}\right) + \mu\left(\frac{n}{3}\right), & n\equiv 0\mod 6, \\
\mu\left(\frac{n}{2}\right), & n\equiv 2,4\mod 6, \\
\mu\left(\frac{n}{3}\right), & n\equiv 3\mod 6. \\
\end{array}
\right.
\end{eqnarray*}
\\

\noindent
Since
\[
-1 = \mu_{\alpha^2}(6) \not= \mu_{\alpha^2}(2) \mu_{\alpha^2}(3) = 0\cdot 0
\]
we have $\mu_{\alpha^2}$ is not multiplicative.
\\

\noindent
We apply $\alpha$--contraction one last time and get
\begin{eqnarray*}
T_{\mu_{\alpha^2},\alpha}(x)
& = &
\sum_{\alpha(n)\leq x} \mu_{\alpha^2}(n) \left\lfloor \frac{x}{\alpha(n)} \right\rfloor \\
& = &
\sum_{\alpha(n)\leq x} (1*\mu_{\alpha^2})(a_n) \\
& = & 
\sum_{\alpha(n)\leq x} \left\lfloor \frac{1}{a_{a_{a_n}}} \right\rfloor \\
& = &
\left\{
\begin{array}{ll}
1, & n < 2, \\
2, & n < 3, \\
3, & n < 4, \\
4, & n \geq 4.
\end{array}
\right.
\end{eqnarray*}
\\

\noindent
With a similar argument as above we get
\[
S_{\mu_{\alpha^2},\alpha}(x) = M(x) + M\left(\frac{x}{2}\right) + M\left(\frac{x}{3}\right)+ M\left(\frac{x}{4}\right)
\]
and therefore
\begin{eqnarray*}
\mu_{\alpha^3}(n)
& = &
\mu(n)
+
\left\{
\begin{array}{ll}
\mu\left(\frac{n}{2}\right), & n\equiv 0\mod 2 \\
0, & n\not\equiv 0\mod 2
\end{array}
\right. 
+
\left\{
\begin{array}{ll}
\mu\left(\frac{n}{3}\right), & n\equiv 0\mod 3 \\
0, & n\not\equiv 0\mod 3
\end{array}
\right.\\
& &
\hspace{1cm} 
+
\left\{
\begin{array}{ll}
\mu\left(\frac{n}{4}\right), & n\equiv 0\mod 4 \\
0, & n\not\equiv 0\mod 4
\end{array}
\right.\\
& = &
\mu(n)
+
\left\{
\begin{array}{ll}
\mu\left(\frac{n}{2}\right) + \mu\left(\frac{n}{3}\right) + \mu\left(\frac{n}{4}\right), & n\equiv 0\mod 12 \\
\mu\left(\frac{n}{2}\right), & n\equiv 2, 10\mod 12 \\
\mu\left(\frac{n}{3}\right), & n\equiv 3, 9\mod 12 \\
\mu\left(\frac{n}{2}\right) + \mu\left(\frac{n}{4}\right), & n\equiv 4, 8\mod 12 \\
\mu\left(\frac{n}{2}\right) + \mu\left(\frac{n}{3}\right), & n\equiv 6\mod 12 \\
\end{array}
\right.\\
\end{eqnarray*}
\begin{eqnarray*}
& = &
\left\{
\begin{array}{ll}
\mu\left(\frac{n}{2}\right) + \mu\left(\frac{n}{3}\right) + \mu\left(\frac{n}{4}\right), & n\equiv 0\mod 12 \\
\mu(n) + \mu\left(\frac{n}{2}\right), & n\equiv 2, 10\mod 12 \\
\mu(n) + \mu\left(\frac{n}{3}\right), & n\equiv 3, 9\mod 12 \\
\mu\left(\frac{n}{2}\right) + \mu\left(\frac{n}{4}\right), & n\equiv 4, 8\mod 12 \\
\mu(n) + \mu\left(\frac{n}{2}\right) + \mu\left(\frac{n}{3}\right), & n\equiv 6\mod 12 \\
\mu(n), & n\equiv 1, 5, 7, 11\mod 12
\end{array}
\right.\\
& = &
\left\{
\begin{array}{ll}
\mu\left(\frac{n}{2}\right) + \mu\left(\frac{n}{4}\right), & n\equiv 0\mod 12 \\
0, & n\equiv 2, 10\mod 12 \\
\mu(n) + \mu\left(\frac{n}{3}\right), & n\equiv 3, 9\mod 12 \\
\mu\left(\frac{n}{2}\right) + \mu\left(\frac{n}{4}\right), & n\equiv 4, 8\mod 12 \\
\mu\left(\frac{n}{3}\right), & n\equiv 6\mod 12 \\
\mu(n), & n\equiv 1, 5, 7, 11\mod 12
\end{array}
\right.\\
& = &
\left\{
\begin{array}{ll}
\mu\left(\frac{n}{2}\right) + \mu\left(\frac{n}{4}\right), & n\equiv 0, 4, 8\mod 12 \\
0, & n\equiv 2, 10\mod 12 \\
\mu(n) + \mu\left(\frac{n}{3}\right), & n\equiv 3, 9\mod 12 \\
\mu\left(\frac{n}{3}\right), & n\equiv 6\mod 12 \\
\mu(n), & n\equiv 1, 5, 7, 11\mod 12
\end{array}
\right.\\
& = &
1, 0, 0, 0, -1, -1, -1, -1, -1, 0, -1, 0, \\
& &
-1, 0, 0, 0, -1, 1, -1, 0, 0, 0, -1, 1, \ldots
\end{eqnarray*}

\noindent
Again we check for multiplicativity and obtain
\[
-1 = \mu_{\alpha^3}(6) \not= \mu_{\alpha^3}(2) \mu_{\alpha^3}(3) = 0\cdot 0
\]
and therefore $\mu_{\alpha^3}$ is not multiplicative.
\\

\noindent
Since $\left\lfloor \frac{1}{a_{a_{a_{a_n}}}} \right\rfloor = \left\lfloor \frac{1}{a_{a_{a_n}}} \right\rfloor$ for $n\in\mathbb{N}$
we have the remarkable identity
\[
\sum_{\alpha(d)=n} \mu_{\alpha^3}(d) =\mu_{\alpha^3}(n)
\]
for $1\leq n\in\mathbb{N}$.
This implies $\mu_{\alpha^m}=\mu_{\alpha^3}$ for $m\geq 3$.
\\

\noindent
The above fixed point property and linearity of $\alpha$--contraction yields
\[
(\mu_{\alpha^3}-\mu_{\alpha^2})_\alpha = \mu_{\alpha^4}-\mu_{\alpha^3} = 0
\]
and thus $\Delta_{23}:= \mu_{\alpha^2} - \mu_{\alpha^3}$ is in the kernel of $\alpha$--contraction.
Let us find an explicit representation for $\Delta_{23}$.
\begin{eqnarray*}
\Delta_{23}(n)
& = &
\mu(n)
+
\left\{
\begin{array}{ll}
\mu\left(\frac{n}{2}\right) + \mu\left(\frac{n}{3}\right), & n\equiv 0\mod 6 \\
\mu\left(\frac{n}{2}\right), & n\equiv 2,4\mod 6 \\
\mu\left(\frac{n}{3}\right), & n\equiv 3\mod 6 \\
\end{array}
\right.\\
& & 
\hspace{1cm}
-
\left\{
\begin{array}{ll}
\mu\left(\frac{n}{2}\right) + \mu\left(\frac{n}{4}\right), & n\equiv 0, 4, 8\mod 12 \\
0, & n\equiv 2, 10\mod 12 \\
\mu(n) + \mu\left(\frac{n}{3}\right), & n\equiv 3, 9\mod 12 \\
\mu\left(\frac{n}{3}\right), & n\equiv 6\mod 12 \\
\mu(n), & n\equiv 1, 5, 7, 11\mod 12
\end{array}
\right.\\
& = &
\mu(n)
+
\left\{
\begin{array}{ll}
\mu\left(\frac{n}{2}\right) + \mu\left(\frac{n}{3}\right), & n\equiv 0,6\mod 12 \\
\mu\left(\frac{n}{2}\right), & n\equiv 2,4,8,10\mod 12 \\
\mu\left(\frac{n}{3}\right), & n\equiv 3,9\mod 12 \\
\end{array}
\right.\\
& &
\hspace{1cm}
-
\left\{
\begin{array}{ll}
\mu\left(\frac{n}{2}\right) + \mu\left(\frac{n}{4}\right), & n\equiv 0, 4, 8\mod 12 \\
0, & n\equiv 2, 10\mod 12 \\
\mu(n) + \mu\left(\frac{n}{3}\right), & n\equiv 3, 9\mod 12 \\
\mu\left(\frac{n}{3}\right), & n\equiv 6\mod 12 \\
\mu(n), & n\equiv 1, 5, 7, 11\mod 12
\end{array}
\right.\\
& = &
\left\{
\begin{array}{ll}
\mu(n) - \mu\left(\frac{n}{3}\right) - \mu\left(\frac{n}{4}\right), & n\equiv 0\mod 12 \\
\mu(n) + \mu\left(\frac{n}{2}\right), & n\equiv 2\mod 12 \\
0, & n\equiv 3, 9\mod 12 \\
\mu(n) - \mu\left(\frac{n}{4}\right), & n\equiv 4, 8\mod 12 \\
\mu(n) + \mu\left(\frac{n}{2}\right), & n\equiv 6\mod 12 \\
\mu(n) + \mu\left(\frac{n}{2}\right), & n\equiv 10\mod 12 \\
0, & n\equiv 1, 5, 7, 11\mod 12
\end{array}
\right.\\
\end{eqnarray*}
\begin{eqnarray*}
& = &
\left\{
\begin{array}{ll}
 - \mu\left(\frac{n}{4}\right), & n\equiv 0\mod 12 \\
0, & n\equiv 2\mod 12 \\
0, & n\equiv 3, 9\mod 12 \\
 - \mu\left(\frac{n}{4}\right), & n\equiv 4, 8\mod 12 \\
0, & n\equiv 6\mod 12 \\
0, & n\equiv 10\mod 12 \\
0, & n\equiv 1, 5, 7, 11\mod 12
\end{array}
\right.\\
& = &
\left\{
\begin{array}{ll}
 - \mu\left(\frac{n}{4}\right), & n\equiv 0,4,8\mod 12 \\
0, & n\equiv 1, 2, 3, 5, 6, 7, 9, 10, 11\mod 12
\end{array}
\right.\\
& = &
\left\{
\begin{array}{ll}
 - \mu\left(\frac{n}{4}\right), & n\equiv 0 \mod 4 \\
0, & n\equiv 1, 2, 3 \mod 4
\end{array}
\right.\\
\end{eqnarray*}

\noindent
Since $f_\alpha(2) = 0$ we know that $\alpha$--contraction as a linear operator on its natural domain, the arithmetic functions, is not surjective.
With $\Delta_{23}$ as a nontrivial kernel element, $\alpha$--contraction is not injective either.

\section{The $\alpha$--contraction of Liouville $\lambda$}

\noindent
In this section we will compute $\lambda_\alpha$. 
To that purpose we need a result of Cohn \cite{Cohn64} and independently by Wyler \cite{Wyler64} ensuring that the only square Fibonacci numbers are $a_1= a_2=1$ and $a_{12}=144$. 
This result was later generalized in various directions cf. eg. London and Finkelstein \cite{LonFin69}, Lagarias and Weissel \cite{LagWei81}, Ribenboim \cite{Rib89} and Bilu, Hanrot and Voutier \cite{BilHanVou01}.
\\

\noindent
As in the examples above we find that $1*\lambda$ is the characteristic function of some set, namely the set of square numbers.
Using the characterisation of square Fibonacci numbers we get 
\begin{eqnarray*}
\lambda_\alpha(n)
& = &
M(n) - M(n-1) + M\left(\frac{x}{2}\right) - M\left(\frac{x-1}{2}\right) + M\left(\frac{x}{12}\right)- M\left(\frac{x-1}{12}\right) \\
& = &
\mu(n)
+
\left\{
\begin{array}{ll}
\mu\left(\frac{n}{2}\right), & n\equiv 0\mod 2 \\
0, & n\not\equiv 0\mod 2
\end{array}
\right.
+
\left\{
\begin{array}{ll}
\mu\left(\frac{n}{12}\right), & n\equiv 0\mod 12 \\
0, & n\not\equiv 0\mod 12
\end{array}
\right.\\
& = &
\mu(n)
+
\left\{
\begin{array}{ll}
\mu\left(\frac{n}{2}\right) + \mu\left(\frac{n}{12}\right), & n\equiv 0\mod 12 \\
\mu\left(\frac{n}{2}\right), & n\equiv 2,4,6,8,10\mod 12 
\end{array}
\right.\\
& = &
\left\{
\begin{array}{ll}
\mu(n), & n \equiv 1,3,5,7,9,11\mod 12\\
\mu\left(\frac{n}{2}\right) + \mu\left(\frac{n}{12}\right), & n\equiv 0\mod 12 \\
\mu(n)+\mu\left(\frac{n}{2}\right), & n\equiv 2,4,6,8,10\mod 12 
\end{array}
\right.
\end{eqnarray*}

\noindent
Observe that 
\[
2=\lambda_\alpha(12) \not= \lambda_\alpha(4)\lambda_\alpha(3) = 1 
\] 
and thus $\lambda_\alpha$ is not multiplicative.

\section{Euler products of $\alpha$--contractions}

\noindent
In general $\alpha$--contractions are not multiplicative and thus Euler products might be hard to find.
However, in the cases considered we can exploit Corollary \arabic{DirichletProdAtFibo} and obtain e.g.
\[
\zeta(s) \sum_{n=1}^\infty \frac{\lambda_\alpha(n)}{n^s} = \sum_{n=1}^\infty \frac{(1*\lambda)(a_n)}{n^s} = 1+ 2^{-s}+12^{-s}
\]
for $s\in\mathbb{C}$ (as analytic extensions).
Therefore
\[
\sum_{n=1}^\infty \frac{\lambda_\alpha(n)}{n^s}
=
(1+ 2^{-s}+12^{-s})
\prod_p (1-p^{-s})
\]
for $s=\sigma +i t$ with real part $\sigma>1$.
With a similar argument we get
\begin{eqnarray*}
\sum_{n=1}^\infty \frac{\mu_\alpha(n)}{n^s} 
& = & (1-4^{-s}) \prod_{p>2} (1-p^{-s}), \\
\sum_{n=1}^\infty \frac{\mu_{\alpha^2}(n)}{n^s} 
& = & (1+2^{-s}+3^{-s}) \prod_{p} (1-p^{-s}), \\
\sum_{n=1}^\infty \frac{\mu_{\alpha^3}(n)}{n^s} 
& = & (1+2^{-s}+3^{-s}+4^{-s}) \prod_{p} (1-p^{-s})
\end{eqnarray*}
again for $s=\sigma +i t$ with real part $\sigma>1$.



\end{document}